\documentclass{article}
\usepackage{amsmath,amsfonts,amsthm,amstext,amssymb,amscd,amsbsy}
\usepackage[latin1]{inputenc}
\usepackage{indentfirst}
\usepackage[dvips]{color}
\usepackage{color}
\usepackage{epsfig}
\usepackage{hyperref}

\newtheorem{teo}{Theorem}
\newtheorem{lema}{Lemma}

\newtheorem{exemplo}{Example}
\newtheorem{defi}{Definition}
\newtheorem{prop}{Proposition}
\newtheorem{corol}{Corollary}
\newcommand{\C}{\mathbb{C}}

\newcommand{\R}{\mathbb{R}}
\newcommand{\dem}{\noindent \textbf{Proof:}\quad}

\newcommand{\gt}{\mathfrak}

\newcommand{\al}{\alpha}

\def\bq{\begin{equation}}
\def\eq{\end{equation}}
\begin{document}

\title{Equigeodesics on full, $G_{2}$ and a rank-three condiction on flag
manifolds}
\author{Neiton P. da Silva \thanks{%
Department of Mathematics, npsilva@ime.unicamp.br}, Lino Grama \thanks{%
Department of Mathematics, linograma@gmail.com} and Caio J.C. Negreiros 
\thanks{%
Department of Mathematics, caione@ime.unicamp.br}. \\
Institute of Mathematics, Statistics and Scientific Computation,\\
P.O. Box 6065, University of Campinas - UNICAMP.}
\date{}
\maketitle

\begin{abstract}
This paper provides a characterization and examples of homogeneous geodesics
on full $G/T$ and $G_{2}$ flag manifolds. We discuss for generalized root
systems the property of sum-zero triple of $T$-roots and give several applications
of this result.
\end{abstract}


\newline



\section{Introduction}

An important class of homogeneous manifolds are the orbits of the adjoint
action of a semisimple compact Lie group, called \emph{generalized flag
manifolds}. Such manifolds can be described by a quotient $\mathbb{F}=G/C(T)$%
, where $C(T)$ is the centralizer of a torus $T$ of the Lie group $G$. If $%
C(T)=T$ then $\mathbb{F}=G/T$ is called \emph{full flag manifold}.

These manifolds were studied by many authors by the 50's, with focus on its
topological properties \cite{Besse}. There are also many recent papers
related to the $G$-invariant geometry in flag manifolds, for instance \cite%
{MN}, \cite{sm}, \cite{arva3}, \cite{CNG} and \cite{Evan2}.

This paper deals with two classical subjects in Riemannian geometry:
geodesics and $G$-invariant geometries.

Let $G/K$ be a homogeneous manifold with origin $o=eK$ (trivial coset) and $%
g $ be a $G$-invariant metric. A geodesic $\gamma(t)$ on $G/K$ through the
origin $o$ is called \textit{homogeneous} if it is the orbit of a $1$%
-parameter subgroup of $G,$ that is, 
\begin{equation*}
\gamma(t)=(\exp tX)\cdot o,
\end{equation*}
where $X \in \mathfrak{g}$. The vector $X$ is called a geodesic vector.

In \cite{CNG} it was introduced the notion of \emph{homogeneous equigeodesics%
}. An \emph{homogeneous equigeodesic} is an homogeneous curve $\gamma $
which is \emph{geodesic with respect to any $G$-invariant metric}. It was
obtained condictions for geometrical flag manifolds (i.e. of type $A_{l})$
to admit homogeneous equigeodesics. All such condictions were described in
terms of equigeodesic vectors.

In this paper we provide a characterization of all homogeneous equigeodesics
in any full flag manifold $G/T$. Such characterization, given in terms of
the equigeodesic vectors (see Section 5 for further details). Our first
result is

\textbf{Theorem A: }\emph{\label{teint1}Let $\mathbb{F}=G/T$ with $T$ a
maximal torus on $G$ and $X\in T_{o}\mathbb{F}$. If $X=X_{\alpha
_{1}}+\ldots +X_{\alpha _{n}}$, with $X_{\alpha _{i}}\in \mathfrak{u}%
_{\alpha _{i}}$ where $\mathfrak{u}_{\alpha _{i}}$ is the root space
associated to the root $\alpha _{i}$ then, the curve $\gamma (t)=(\exp
tX)\cdot o$ is an equigeodesic if, and only if, $\alpha _{i}\pm \alpha _{j}$
are not roots, for any $i,j\in \{1,\ldots ,n\}.$ } \newline

After Theorem A\textbf{,} we studied homogeneous equigeodesics in flag
manifolds of type $G_{2}$ obtaining its classification. According to \cite{BFR86} these manifolds are classified in three types: $%
G_{2}(\alpha _{1})=G_{2}/U(2)$, where $U(2)$ is represented by the long
root; $G_{2}(\alpha _{2})=G_{2}/U(2)$, where $U(2)$ is represented by the
short root; and $G_{2}/T$, where $T$ is a maximal torus of $G_{2}$.
According Theorem A, we only have to discuss the cases $G_{2}(\alpha _{1})$
and $G_{2}(\alpha _{2})$.

The tangent space at origin $o$ (trivial coset) of $G_{2}(\alpha _{2})$ can
be written as 
\begin{equation*}
T_{o}G_{2}(\alpha _{2})=\mathfrak{m}_{\alpha _{1}}\oplus \mathfrak{m}%
_{2\alpha _{1}}
\end{equation*}%
where $\mathfrak{m}_{\alpha _{1}},\mathfrak{m}_{2\alpha _{1}}$ are the
irreducible submodules of the isotropic representation. Thus the next two
results classify all the\emph{\ equigeodesic}\textit{s} on $G_{2}(\alpha
_{1})$ and $G_{2}(\alpha _{2}).$

\textbf{Theorem B: }\emph{Let $G_{2}(\alpha _{2})=G_{2}/U(2)$. Then $X$ is
an \emph{equigeodesic vector} if, and only if, $X\in \mathfrak{m}_{\alpha
_{1}}$ or $X\in \mathfrak{m}_{2\alpha _{1}}$. }\newline

On the other hand for the flag manifolds $G_{2}(\alpha _{1})$ we prove the
following

\textbf{Theorem C: }\emph{A vector $X\in T_{o}G_{2}(\alpha _{1})$ is an
equigeodesic vector iff the coefficients of $X$ are solutions of a
non-linear algebraic system of equations. Such system of equations can be
solved explicitly. }\newline

In section 5 we also provide several examples of homogeneous equigeodesics
in any class of full and $G_{2}$ type of flag manifolds: $G/T,$ $%
G_{2}(\alpha _{1})$, $G_{2}(\alpha _{2})$.

One key point in the understanding of the invariant Hermitian geometry of
flags is the study of the behavior of triples of roots (the relevance of
this fact it was first noticed in \cite{MN}). The concept of sum-zero triple
for root systems it was introduced in \cite{ncsa} and is naturally
associated to the study of (1,2)-symplectic metrics on flags.

In this paper following \cite{arva3} we generalize the notion
of root systems for an arbitrary flag manifold $\mathbb{F}=G/K$ and call it
by a system of T-roots $R_{T}$. It is defined as the restriction of the root
system $R$ of the corresponding Lie algebra $\mathfrak{g}$ to the center $%
\mathfrak{t}$ of the (stability) subalgebra $\mathfrak{k}$ of $K$.

The following result is useful in order to determine the set of $T$-roots 
\emph{\ }and is connected to properties of Einstein metrics on flag
manifolds.


\textbf{Theorem D: }\emph{\label{teint2}Let $M=G/K$ and \ $R_{T}$ a
corresponding set of $T$-roots. If $R_{T}$ contains more than one positive $%
T $-root then every $T$-root belongs to some $T$-zero sum triple. }\newline

In a forthcoming paper \cite{neitonnir} we will apply this result in order
to obtain a description of Hermitian classes in terms of $T$-roots. In $\cite%
{sm}$ all the invariant Hermitian structures were classified on full flag
manifolds. Among all these metrics it is natural to determine the ones that
are Einstein.

We can related a result in (\cite{Wang e Ziller}, Corollary 1.5) concerning
to the normal metric with theorem D and a result we have derived connecting
Einstein metrics with the dimensions of the modules $\mathfrak{m}_{i}.$

\textbf{Corollary E: }\emph{Let $\mathbb{F}=G/K$ and $T_{o}M=\mathfrak{m}%
_{1}\oplus\mathfrak{m}_{2}$ then $R_{T}=\left\{ \pm\zeta ,\pm 2\zeta\right\} 
$ where $\zeta\in$ $R_{T}-\left\{ 0\right\} .$ Furthermore, if $\dim$ $%
\mathfrak{m}_{1}\neq\dim$ $\mathfrak{m}_{2}$ then any invariant Eintein
metric on $\mathbb{F}$ satisfies $\lambda_{1}\neq\lambda _{2}$. }\newline

We, just for completeness, derive the Einstein equations for $G_{2}/T$ and
describe explicitly the invariant Kähler-Einstein metric correponding to
each invariant complex structures on $G_{2}/T$. \newline

The paper is organized in the following format. In the Sections 2 and 3 we
summarize some results about the geometry of the flag manifolds, describe
the isotropy representation, $T$-roots, invariant metrics and the invariant
almost-complex structures. In the Section 4 we present the classification of
the flag manifolds of type $G_{2}$. In the Section 5 we prove a characterization of equigeodesic vector on full and $G_2$ flag manifolds and give several examples of such curves. Finally, in Section 6, we present the result about the zero sum triple of $T$-roots, with
applications in the study of the invariant Einstein metrics.

\section{Flag manifolds}

In this section we briefly review some basic facts on the structure of
homogeneous spaces, flag manifolds and describe the T-root system. 
\newline

\emph{I. Homogeneous spaces.} Consider the homogeneous manifold $M=G/K$ with 
$G$ a compact semi-simple Lie group and $K$ a closed subgroup. Let $%
\mathfrak{g}$ and $\mathfrak{k}$ be the corresponding Lie algebras. The
Cartan-Killing form $\left\langle ,\right\rangle$ is nondegenerate and
negative definite in $\mathfrak{g}$, thus giving rise to the direct sum
decomposition $\mathfrak{g=k\oplus m}$ where $\mathfrak{m}$ is $\mbox{Ad}(K)$%
-invariant. 
We may identify $\mathfrak{m}$ with the tangent space $T_oM$ at $o=eK$. The
isotropy representation of a reductive homogeneous space is equivalent to
the homomorphism $j:K\rightarrow GL(T_oM)$ given by $j(k)=\mbox{Ad}(k)\big|_{%
\mathfrak{m}}$.



\emph{II. Generalized flag manifolds.} A homogeneous space $\mathbb{F}=G/K$
is called a generalized flag manifold if $G$ is simple Lie group and the
isotropy group $K$ is the centralizer of a one-parametrer subgroup of $G$, $%
\exp tw$ $(w\in \mathfrak{g})$. Equivalently, $\mathbb{F}$ is an adjoint
orbit $\mbox{Ad}(G)w$, where $w\in \mathfrak{g}$. The generalized flag
manifolds (also refereed to as a Kählerian $C$-spaces) have been classified
in \cite{BFR86},\cite{Wa}.

Here the direct sum decomposition $\mathfrak{g}=\mathfrak{k}\oplus\mathfrak{m%
}$ has a more complete description. Let $%
\mathfrak{h}^{\mathbb{C}}$ be a Cartan subalgebra of the complexification $%
\mathfrak{k}^{\mathbb{C}}$ of $\mathfrak{k}$, which is also a Cartan
subalgebra of $\mathfrak{g}^{\mathbb{C}}$. Let $R$ and $R_K$ be the root
systems of $\mathfrak{g}^{\mathbb{C}}$ and $\mathfrak{k}^{\mathbb{C}}$,
respectively, and $R_M=R\backslash R_K$ be the set of complementary roots.
We have the Cartan decompositions 
\begin{equation*}
\mathfrak{g}^{\mathbb{C}}=\mathfrak{h}^{\mathbb{C}} \oplus \sum_{\alpha\in R}%
{\mathfrak{g}_{\alpha}},\hspace{1cm} \mathfrak{k}^{\mathbb{C}}=\mathfrak{h}^{%
\mathbb{C}} \oplus \sum_{\alpha\in R_K}{\mathfrak{g}_{\alpha}}, \hspace{1cm} 
\mathfrak{m}^{\mathbb{C}}=\sum_{\alpha\in R_M}{\mathfrak{g}_{\alpha}}
\end{equation*}
where $\mathfrak{g}_{\alpha}$ denotes the root space associated to root $%
\alpha$ and $\mathfrak{m}^{\mathbb{C}}$ is isomorphic to $(T_o\mathbb{F})^{%
\mathbb{C}}$ and $\mathfrak{h}=\mathfrak{h}^{\mathbb{C}}\cap\mathfrak{g}$.

We fix a Weyl basis in $\mathfrak{m}^\mathbb{C}$, namely, elements $%
E_{\alpha}\in \mathfrak{g}_{\alpha}$ such that $\left\langle
E_{\alpha},E_{-\alpha}\right\rangle=1$ and $[E_{\alpha},E_{\beta}]=m_{%
\alpha,\beta}E_{\alpha+\beta}$, with $m_{\alpha,\beta}\in \mathbb{R}$, $%
m_{\alpha,\beta}=-m_{\beta,\alpha}$, $m_{\alpha,\beta}=-m_{-\alpha,-\beta}$
and $m_{\alpha,\beta}=0$ if, and only if, $\alpha+\beta$ is not a root. The
corresponding \emph{real} Weyl basis in $\mathfrak{m}$ consists of the
vectors $A_{\alpha}=E_{\alpha}-E_{-\alpha}$, $S_{\alpha}=E_{\alpha}+E_{-%
\alpha}$ and $\mathfrak{u}_{\alpha}=\mbox{span}_{\mathbb{R}}\{
A_{\alpha},iS_{\alpha}\}$, where $\alpha\in R^+$, the set of positive roots.%
\newline

The real tangent space $T_o\mathbb{F}$ is naturally identified with 
\begin{equation*}
\displaystyle{\mathfrak{m}=\bigoplus_{\alpha\in R_M^+} \mathfrak{u}_{\alpha}}%
.
\end{equation*}

Some of the spaces $\mathfrak{u}_\alpha$ are not $\mbox{Ad}(K)$-modules,
unless $\mathbb{F}$ is a full flag manifold. To get the \emph{irreducible} $%
Ad(K)$-modules, we proceed as in \cite{ale1}. Let 
\begin{equation*}
\mathfrak{t}=Z(\mathfrak{k}^{\mathbb{C}})\cap\mathfrak{h}=\{ X\in\mathfrak{h}%
:\phi(x)=0 \,\, \forall \phi \in R_K\}.
\end{equation*}
If $\mathfrak{h}^*$ and $\mathfrak{t}^*$ are the dual space of $\mathfrak{h}$
and $\mathfrak{t}$ respectively, we consider the restriction map 
\begin{equation}  \label{projkappa}
\kappa:\mathfrak{h}^*\rightarrow \mathfrak{t}^*, \quad\quad\quad\quad
\kappa(\alpha)=\alpha|_{\mathfrak{t}}
\end{equation}
and set $R_T=\kappa(R_M)$. The elements on $R_T$ are called $T$-roots. The
irreducible $\mbox{ad}(\mathfrak{k}^{\mathbb{C}})$-invariant sub-modules of $%
\mathfrak{m}^{\mathbb{C}},$ and the corresponding irreducible sub-modules
for the $\mbox{ad}(\mathfrak{k})$-module $\mathfrak{m},$ are given by 
\begin{equation*}
\mathfrak{m}^\mathbb{C}_{\xi}=\sum_{\kappa(\alpha)=\xi} {\mathfrak{g}%
_{\alpha}}\quad\quad(\xi\in R_T),\quad\quad\quad\quad \mathfrak{m}_{\eta}=%
\displaystyle{\sum_{\kappa(\alpha)=\eta}{\mathfrak{u}_{\alpha}}}
\quad\quad(\eta\in R^+_T=\kappa(R_M^+)).
\end{equation*}
Hence we have the direct sum of complex and real irreducible modules, 
\begin{equation*}
\mathfrak{m}^\mathbb{C}=\displaystyle{\sum_{\eta\in R_T}{\mathfrak{m}^%
\mathbb{C}_{\eta}}},\quad\quad\quad \mathfrak{m}=\displaystyle{\sum_{\eta\in
R^+_T}{\mathfrak{m}_{\eta}}}.
\end{equation*}

\section{Invariant metrics and iacs}

An invariant metric $g$ on $\mathbb{F}$ is uniquely determined by a scalar
product $B$ on $\mathfrak{m}$ of the form 
\begin{equation*}
B(\cdot\,,\cdot)=-\left\langle \Lambda \cdot\,,\cdot \right\rangle=
\lambda_1(-\left\langle \cdot\,,\cdot \right\rangle)|_{\mathfrak{m}%
_1}+\ldots + \lambda_s(-\left\langle \cdot\,,\cdot \right\rangle)|_{%
\mathfrak{m}_s,}
\end{equation*}
where the linear map $\Lambda:\mathfrak{m}\rightarrow \mathfrak{m}$ is
symmetric, positive-definite with respect to the Cartan-Killing form, $%
\lambda_i>0$ and $\mathfrak{m}_i$ are the irreducible $\mbox{Ad}(K)$%
-sub-modules. Each $\mathfrak{m}_i$ is an eigenspace of $\Lambda$
corresponding to the eigenvalue $\lambda_i$. In particular, the vectors $%
A_{\alpha},S_{\alpha}$ of the real Weyl basis are eigenvectors of $\Lambda$
correponding to the same eigenvalue $\lambda_{\alpha}$. We abuse of notation
and say that $\Lambda$ itself is an invariant metric.

The inner product $B$ admits a natural extension to a symmetric bilinear
form on the complexification $\mathfrak{m}^\mathbb{C}$ of $\mathfrak{m}$. We
do not change notation for these objects in $\mathfrak{m}$ and $\mathfrak{m}%
^{\mathbb{C}}$ either for the bilinear form $B$ or for the corresponding
complexified map $\Lambda$.

It is well known that an $G$-invariant almost complex structure (abbreviated
iacs) on $\mathbb{F}$ is completely determined by its value $J\colon%
\mathfrak{m}\rightarrow\mathfrak{m}$ in the tangent space at the origin. The
linear endomorphism $J$ satisfies $J^2=-1$ and $Ad(K)J=JAd(K)$. We will also
denote by $J$ its complexification to $\mathfrak{m}^{\mathbb{C}}$. The
eigenvalues of $J$ are $\pm i$ and the corresponding eigenvectors are
denoted by $T^{(1,0)}_{o}\mathbb{F}=\{X\in T_{o}\mathbb{F}:JX=iX\}$ and $%
T^{(0,1)}_{o}\mathbb{F}=\{X\in T_{o}\mathbb{F}:JX=-iX\}$. Thus we have the
following decomposition of the tangent complex space at the origin $%
\mathfrak{m}^{\mathbb{C}}=T^{(1,0)}_{o}\mathbb{F}\oplus T^{(0,1)}_{o}\mathbb{%
F}$. The invariance of $J$ entails that $J(\mathfrak{g}^{\mathbb{C}%
}_{\alpha})=\mathfrak{g}^{\mathbb{C}}_{\alpha}$ for all $\alpha\in R$. Then $%
JE_{\alpha}=i\varepsilon_{\alpha}E_{\alpha}$, with $\varepsilon_{\alpha}=\pm
1$.

As $A_{\alpha}=E_{\alpha}-E_{-\alpha}$ and $S_{\alpha}=E_{\alpha}+E_{-%
\alpha} $ we obtain that $iA_{\alpha}=-iA_{-\alpha}$ and $%
S_{\alpha}=S_{-\alpha}$. Note that $E_{-\alpha}=\frac{1}{2}%
(i(iA_{\alpha})+S_{\alpha})$, then 
\begin{equation*}
i\varepsilon_{-\alpha}\frac{1}{2}(i(iA_{\alpha})+S_{\alpha})
=i\varepsilon_{-\alpha}E_{-\alpha}=JE_{-\alpha}=\frac{1}{2}%
(iJ(iA_{\alpha})+J(S_{\alpha})).
\end{equation*}
Comparing the reals and imaginary terms in the left side of the first
equation and in the right side of the third equation, we obtain $%
J(iA_{\alpha})=\varepsilon_{-\alpha}S_{\alpha}$ and $J(S_{\alpha})=-%
\varepsilon_{-\alpha}(iA_{\alpha})$, then 
\begin{equation*}
\varepsilon_{-\alpha}S_{\alpha}=J(iA_{\alpha})=-J(iA_{-\alpha})=-%
\varepsilon_{\alpha}S_{-\alpha},
\end{equation*}
so $\varepsilon_{\alpha}=-\varepsilon_{-\alpha}$, with $\alpha\in R$.

The irreducible $ad(\mathfrak{k}^{\mathbb{C}})$-modules are $\mathfrak{m}^{%
\mathbb{C}}_i$ invariant by $J$, that is, $J\mathfrak{m}^{\mathbb{C}}_i=%
\mathfrak{m}^{\mathbb{C}}_i$. Then using the Schur's Lemma we get 
\begin{equation*}
J=i\varepsilon_{1}Id|_{\mathfrak{m}^{\mathbb{C}}_1}\oplus\cdots\oplus
i\varepsilon_{r}Id|_{\mathfrak{m}^{\mathbb{C}}_r}.
\end{equation*}
Moreover, if $\delta$ is a $T$-root we have $\varepsilon_{\delta}=%
\varepsilon_{\alpha}=-\varepsilon_{-\alpha}=-\varepsilon_{-\delta}$ where $%
\alpha$ is any root in $R$ such that $k(\alpha)=\delta$. Thus we obtain

\begin{prop}
Let $\mathbb{F}$ be a flag manifold and $R_{T}$ the correspondent set of
T-roots. Then any iacs $J$ on $\mathbb{F}$ is completely determined by a set
of sign $\{\varepsilon_{\delta}, \delta\in R_{T}\}$ ($\varepsilon_{\delta}=%
\pm 1$) satisfying $\varepsilon_{\delta}=-\varepsilon_{-\delta}$ with $%
\delta\in R_{T}$. In particular, $J$ is determined by exactly $|R_{T}^+|$
signs.
\end{prop}

An iacs $J$ is integrable if, and only if, it is torsion free, that is, 
\begin{equation*}
[JX,JY]=[X,Y]+J[X,JY]+J[JX,Y]\hspace{0.5cm}X,Y\in\mathfrak{m}, 
\mbox{   
(see, for instance, \cite{KN})}.
\end{equation*}

\section{Flag manifolds of $G_2$ type.}

We recall some basics facts about the Lie algebra of $G_2$. We can realize
the Lie algebra of $G_2$ as the Lie algebra $\mathfrak{g=sl}\mathbb{(}3,%
\mathbb{C)}\mathbb{\oplus C}^{3}\oplus\left(\mathbb{C}^{3}\right) ^{\ast}$.
A Cartan subalgebra $\mathfrak{h}$ of diagonal matrices on $\mathfrak{sl}%
\mathbb{(}3,\mathbb{C)}$ is also a Cartan subalgebra on $\mathfrak{g}$.

Consider the linear functional $\varepsilon _{i}$ de $\mathfrak{h}$ defined
by: $\varepsilon _{i}\colon diag\{a_{1},a_{2},a_{3}\}\mapsto a_{i}$. A basis
for the root system relative to $(\mathfrak{g},\mathfrak{h})$ is given by $%
\Sigma =\{\alpha _{1}=\varepsilon _{1}-\varepsilon _{2},\alpha
_{2}=\varepsilon _{2}\}$. The corresponding positive roots are $%
R^{+}=\{\alpha _{1}=\varepsilon _{1}-\varepsilon _{2},\text{ }\alpha
_{2}=\varepsilon _{2},\text{ }\alpha _{1}+\alpha _{2}=\varepsilon _{1},\text{
}\alpha _{1}+2\alpha _{2}=-\varepsilon _{3},\text{ }\alpha _{1}+3\alpha
_{2}=\varepsilon _{2}-\varepsilon _{3},\text{ }2\alpha _{1}+3\alpha
_{2}=\varepsilon _{1}-\varepsilon _{3}\}$. The Cartan-Killing form $(\cdot
,\cdot )$ on $\mathfrak{h}^{\ast }$ is given by: 
\begin{align*}
(2\alpha _{1}+3\alpha _{2},2\alpha _{1}+3\alpha _{2})& =(\alpha _{1}+3\alpha
_{2},\alpha _{1}+3\alpha _{2})=(\alpha _{1},\alpha _{1})=\frac{1}{4}\,\,\,%
\text{ (long roots) } \\
(\alpha _{1}+2\alpha _{2},\alpha _{1}+2\alpha _{2})& =(\alpha _{1}+\alpha
_{2},\alpha _{1}+\alpha _{2})=(\alpha _{2},\alpha _{2})=\frac{1}{12}\,\,\,%
\text{ (short roots) }.
\end{align*}

According to \cite{BFR86} there are only three
non-equivalent classes of $G_{2}$ flag manifolds.

The following table list these manifolds.

\begin{center}
\begin{tabular}{|l|c|}
\hline
Flag manifold & $\mathfrak{m}=\oplus_{i=1}^t \mathfrak{m}_i$ \\ \hline\hline
$G_2(\alpha_1)= G_2/U(2)$, where $U(2)$ is represented by the long root & $%
t=3$ \\ \hline
$G_2(\alpha_2)= G_2/U(2)$, where $U(2)$ is represented by the short root & $%
t=2$ \\ \hline
$G_2/T$, where $T$ is a maximal torus of $G_2$ & $t=6$ \\ \hline
\end{tabular}
\end{center}

The second column represents the number of the irreducible non-equivalent
submodules of the isotropic representation.

\section{Homogeneous geodesics and equigeodesics}

In this section we give a characterization of homogeneous equigeodesics in
full flag manifolds and flag manifolds of $G_2$. We start with an
definition.

\begin{defi}
\label{defihomo} Let $(M=G/K,g)$ be a homogeneous Riemannian manifold. A
geodesic $\gamma(t)$ on $M$ through the origin $o$ is called \textit{%
homogeneous} if it is the orbit of a $1$-parameter subgroup of $G,$ that is, 
\begin{equation*}
\gamma(t)=(\exp tX)\cdot o,
\end{equation*}
where $X \in \mathfrak{g}$. The vector $X$ is called a geodesic vector.
\end{defi}

A useful result of Kowalski and Vanhecke \cite{K2} gives an algebraic
characterization of the geodesic vectors.

\begin{teo}
\label{hecke}If $g$ is a $G$-invariant metric, a vector $X \in \mathfrak{g}%
\setminus\{ 0 \}$ is a geodesic vector if, and only if, 
\begin{equation}  \label{eqnfund}
g(X_{\mathfrak{m}},[X,Z]_{\mathfrak{m}})=0,
\end{equation}
for all $Z\in \mathfrak{m}$.
\end{teo}

An important class of the homogeneous manifolds are the \emph{g.o. manifolds}
(geodesic orbit manifold). We say that a homogeneous manifolds is a g.o
manifold if it admits a invariant Riemannian metric such that \emph{all}
geodesics are \emph{homogeneous}. Examples of these manifold are the
homogeneous space equipped with the normal metric and the symmetric spaces.

It is well know that neither full flag manifolds $G/T$ nor flag manifolds of 
$G_2$ admits a left invariant metric (not homotetic to normal metric) such
that all geodesics are homogeneous (see \cite{arva3}).

On other hand, every homogeneous manifold admits at least one homogeneous
geodesics, see \cite{K1}. In case of the group $G$ is semi-simple we have
the following existence result:

\begin{teo}[\protect\cite{K1}]
\label{givevector} If $G$ is semi-simple then $M=G/K$ admits at least $%
m=dim(M)$ mutually orthogonal homogeneous geodesics through the origin $o$.
\end{teo}

An \emph{homogeneous equigeodesic} is an homogeneous curve that is geodesic
with respect to any invariant metric. In \cite{CNG} it was proved that any
flag manifold admits equigeodesics. The following algebraic characterization
of equigeodesics is given in terms of the \emph{equigeodesic vectors}, that
is, vectors $X \in \mathfrak{m}$ such that the orbit $\gamma(t)=(\exp
tX)\cdot o$ is an homogeneous equigeodesic.

\begin{prop}[\protect\cite{CNG}]
\label{proplegal} Let $\mathbb{F}$ be a flag manifold, with $\mathfrak{m}$
isomorphic to $T_o F$. A vector $X\in \mathfrak{m}$ is equigeodesic vector
if, and only if, 
\begin{equation}  \label{eqnlegal}
[X,\Lambda X]_{\mathfrak{m}}=0,
\end{equation}
for any invariant metric $\Lambda$.
\end{prop}

\noindent \textbf{Proof:}\quad Let $g$ be the metric associated with $%
\Lambda $. For $X,Y\in \mathfrak{m}$ we have 
\begin{eqnarray*}
g(X,[X,Y]_{\mathfrak{m}}) &=& -\left\langle \Lambda X,[X,Y]_{\mathfrak{m}}
\right\rangle = - \left\langle \Lambda X,[X,Y] \right\rangle = -\left\langle
[X,\Lambda X],Y \right\rangle,
\end{eqnarray*}
since the decomposition $\mathfrak{g=m+h}$ is $\left\langle ,\right\rangle$%
-orthogonal and the Killing form is $\mbox{Ad}(G)$-invariant, i.e., $%
\mbox{ad}(X)$ is skew-Hermitian with respect to $<,>$. Therefore $X$ is
equigeodesic iff $[X,\Lambda X]_{\mathfrak{m}}=0$ for any invariant scalar
product $\Lambda$.\qed
\newline

We will now give a full characterization of the equigeodesics or
equivalently equigeodesic vectors on \emph{any} full flag manifold $G/T$,
where $G$ is a compact, connected and simple Lie group and $T$ is a maximal
torus on $G$. We recall that the irreducible submodules of the isotropy
representation in full flag manifolds coincides with $\mathfrak{u}%
_\alpha=span_{\mathbb{R}}\{A_\alpha,iS_\alpha \}, \alpha\in R^+$.

\begin{teo}
\label{geodesicfull} Let $G$ be a compact, connected and simple Lie group
with Lie algebra $\mathfrak{g}$, $T$ a maximal torus in $G$ and $G/T$ the
corresponding full flag manifold. Let $X\in \mathfrak{u}_\alpha$, $Y\in 
\mathfrak{u}_\beta$ be nonzero vectors. Then $X+Y\in \mathfrak{m}$ is an
equigeodesic vectors if, and only if, $\alpha \pm \beta$ are not roots.
\end{teo}

\noindent \textbf{Proof:}\quad Let $\{ E_{\alpha}\}_{\alpha\in R}$ be the
Weyl´s basis of $\mathfrak{g}^{\mathbb{C}}$ the complexification of the real
simple algebra $\mathfrak{g}$, and set $X=a_1A_\alpha +b_1iS_\alpha \in 
\mathfrak{u}_\alpha, Y=a_2A_\beta +b_2iS_\beta \in \mathfrak{u}_\beta $,
with $a_1,a_2,b_1,b_2 \in \mathbb{R}$. Using the Weyl´s basis we can write $%
X=cE_{\alpha}-\overline{c}E_{-\alpha}$ and $Y=dE_{\beta}-\overline{d}%
E_{-\beta}$, where $c=a_1+ib_1$ and $d=a_2+ib_2$. Then, 
\begin{eqnarray}
[X+Y,\Lambda(X+Y)]_{\mathfrak{m}}&=&[cE_{\alpha}-\overline{c}%
E_{-\alpha}+dE_{\beta}-\overline{d}E_{-\beta}, \,
c\lambda_{\alpha}E_{\alpha}-\overline{c}\lambda_{\alpha}E_{-\alpha}+d%
\lambda_{\beta}E_{\beta}-\overline{d}\lambda_{\beta}E_{-\beta} ]  \notag \\
&=& (\lambda_{\beta} - \lambda_{\alpha})cd\,
m_{\alpha,\beta}E_{\alpha+\beta} - (\lambda_\beta-\lambda_\alpha) c 
\overline{d}\,m_{\alpha,-\beta}E_{\alpha-\beta}  \notag \\
&& - (\lambda_\beta-\lambda_\alpha) \overline{c}d\,
m_{-\alpha,\beta}E_{-\alpha+\beta} +(\lambda_\beta - \lambda_\alpha)%
\overline{c}\overline{d}\lambda_{\beta}\,m_{-\alpha,-\beta}E_{-\alpha-\beta}.
\label{colchetefinal}
\end{eqnarray}

Suppose that $X+Y$ is equigeodesic. Then $[X+Y,\Lambda (X+Y)]_{\mathfrak{m}%
}=0$ for any invariant metric $\Lambda $ and from equation %
\eqref{colchetefinal} we have $m_{\alpha ,\beta }=-m_{-\alpha ,-\beta }=0$ e 
$m_{\alpha ,-\beta }=-m_{-\alpha ,\beta }=0$ because $c$ and $d$ are nonzero
and therefore $\alpha \pm \beta $ are not roots.

On other hand, suppose that $\alpha \pm \beta $ are not roots. Then $%
m_{\alpha ,\beta }=-m_{-\alpha ,-\beta }=0$ and $m_{\alpha ,-\beta
}=-m_{-\alpha ,\beta }=0$ and from \eqref{colchetefinal} we have $%
[X+Y,\Lambda (X+Y)]_{\mathfrak{m}}=0$ for any invariant metric $\Lambda $
and $X+Y$ is an equigeodesic vector. \qed

\begin{corol}
With the hypothesis from Theorem \ref{geodesicfull}, let $%
X=X_{\alpha_1}+\ldots +X_{\alpha_r}$ such that $X_{\alpha_i}\in \mathfrak{u}%
_{\alpha_i}$ for all $i$. Then $X$ is an equigeodesic if, and only if, $%
\alpha_{p}\pm\alpha_{q}$ are not roots for every $p,q\in \{1,\ldots,r \}$.
\end{corol}

\noindent \textbf{Proof:}\quad Apply Theorem above and the linearity of the
Lie bracket. \qed

\begin{exemplo}[\cite{CNG}]Consider the Lie algebra $A_l=\gt{sl}(l,\C)$. The Cartan sub-algebra of $\gt{sl}(n,\C)$ can be identified with $\gt{h}=\{\mbox{diag}(\varepsilon_1,\ldots,\varepsilon_n);\varepsilon_i\in\C,\sum{\varepsilon_i}=0\}$. 

The root system of the Lie algebra of $\gt{sl}(n)$ has the form $R=\{\alpha_{ij}=\varepsilon_i-\varepsilon_j:i\neq j\}$ and the subset of positive roots is  $R^+=\{\alpha_{ij}:i<j\}$, see \cite{H}. Therefore, $\al_{ij}\pm \al_{pq}$ is not a root if, and only if, $i,j,p,q$ are all distinct. 

\end{exemplo}

\begin{exemplo}
Consider the Lie algebra $\gt{g}^{\C}$ over $\C$, $R$ beeing an associated root system, and $\Sigma$ a simple root system. Two simple roots are said to be {\em orthogonal} if they are not joined in the Dynkin diagram. If $\al_1$ and $\al_2$ are two orthogonal simple roots then $\al_1\pm\al_2$ are not roots, see \cite{H}. 

For example, on the full flag manifold $SO(16)/T$ we have $\gt{g}^{\C}=\gt{so}(16,\C)$ (a Lie algebra of type $D_l$) and the associated Dynkin diagram is given by  
\\

%
%
%
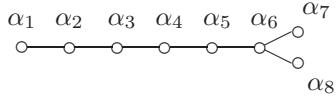
\begin{picture}(160,40)(-15,-23)
\put(0, 0){\circle{4}}
\put(0,10){\makebox(0,0){$\al_1$}}
\put(2, 0){\line(1,0){14}}
\put(18, 0){\circle{4}}
\put(18,10){\makebox(0,0){$\al_2$}}
\put(20, 0){\line(1,0){14}}
\put(36, 0){\circle{4}}
\put(38,10){\makebox(0,0){$\al_3$}}
\put(38, 0){\line(1,0){14}}

\put(54, 0){\circle{4}}
\put(56,10){\makebox(0,0){$\al_4$}}
\put(56, 0){\line(1,0){14}}

\put(72, 0){\circle{4}}
\put(74,10){\makebox(0,0){$\al_5$}}
\put(74, 0){\line(1,0){14}}

\put(90, 0){\circle{4}}
\put(92,10){\makebox(0,0){$\al_6$}}

\put(92, 1){\line(2,1){10}}
\put(92, -1){\line(2,-1){10}}
\put(104.5, 6){\circle{4}}
\put(104.5, -6){\circle{4}}
\put(111.5, 14){\makebox(0,0){$\al_7$}}
\put(108, -16){$\al_8$}

\end{picture}

Hence, any element in the set $\gt{u}_{\al_1}\oplus\gt{u}_{\al_3}$ is an equigeodesic vector since $\al_1\pm\al_3$ are not roots. In the same way, any element in the set $\gt{u}_{\al_2}\oplus\gt{u}_{\al_4}\oplus \gt{u}_{\al_7}$ is equigeodesic vector.

In a completely similar way we find equigeodesic vectors in any full flag manifold.
\end{exemplo}

\begin{exemplo}
Let $G_2/T$ be the full flag manifold and consider the roots $%
\alpha=\alpha_2 $ and $\beta=2\alpha_1+3\alpha_2$. Then any element in the
set $\mathfrak{u}_{\alpha}\oplus\mathfrak{u}_{\beta}$ is an equigeodesic
vector.
\end{exemplo}

We will now give a complete description of the equigeodesic vectors in flag
manifolds of $G_2$.

\emph{a) Let $G_{2}(\alpha _{2})=G_{2}/U(2)$ represented by the short root.}
In this case we have $R_{K}=\{\alpha _{2}\}$ and $R_{M}=\{\alpha _{1}\}$.
The isotropic representation have two irreducible submodules and the tangent
space at $o$ splits as $\mathfrak{m}=\mathfrak{m}_{\alpha _{1}}\oplus 
\mathfrak{m}_{2\alpha _{1}}$, where $\mathfrak{m}_{\alpha _{1}}=\mathfrak{u}%
_{\alpha _{1}}\oplus \mathfrak{u}_{\alpha _{1}+\alpha _{2}}\oplus \mathfrak{u%
}_{\alpha _{1}+2\alpha _{2}}\oplus \mathfrak{u}_{\alpha _{1}+3\alpha _{2}}$
and $\mathfrak{m}_{2\alpha _{1}}=\mathfrak{u}_{2\alpha _{1}+3\alpha _{2}}$.

\begin{lema}
\label{inclusao dois somandos} Let $\mathfrak{m}_{\alpha_1},\mathfrak{m}%
_{\alpha_2}$ be the irreducibles submodules of the isotropy representation
of the flag manifold $G_2(\alpha_2)$. Then $[\mathfrak{m}_{\alpha_1},%
\mathfrak{m}_{2\alpha_1}]\subset \mathfrak{m}_{\alpha_1}$.
\end{lema}

\noindent \textbf{Proof:}\quad Let $\{E_{\alpha }\}_{\alpha \in R_{M}}$ be a
Weyl's basis of $\mathfrak{g_{2}}$ and let $X\in \mathfrak{m}_{\alpha _{1}}$
and $Y\in \mathfrak{m}_{2\alpha _{1}}$. Writing 
\begin{equation*}
X=a_{1}E_{\alpha _{1}}+a_{2}E_{\alpha _{1}+\alpha _{2}}+a_{3}E_{\alpha
_{1}+2\alpha _{2}}+a_{4}E_{\alpha _{1}+3\alpha _{2}}+b_{1}E_{-\alpha
_{1}}+b_{2}E_{-(\alpha _{1}+\alpha _{2})}+b_{3}E_{-(\alpha _{1}+2\alpha
_{2})}+b_{4}E_{-(\alpha _{1}+3\alpha _{2})},
\end{equation*}%
\begin{equation*}
Y=c_{1}E_{2\alpha _{1}+3\alpha _{2}}+c_{2}E_{-(2\alpha _{1}+3\alpha _{2})},
\end{equation*}%
with $b_{i}=-\overline{a_{i}}$ e $c_{2}=-\overline{c_{1}}$, we have 
\begin{eqnarray}
\lbrack X,Y] &=&a_{1}c_{1}E_{-(\alpha _{1}+3\alpha
_{2})}+a_{2}c_{2}E_{-(\alpha _{1}+2\alpha _{2})}+a_{3}c_{2}E_{-(\alpha
_{1}+\alpha _{2})}+a_{4}c_{2}E_{-\alpha _{1}}  \notag \\
&&-b_{1}c_{1}E_{\alpha _{1}+3\alpha _{2}}-b_{2}c_{1}E_{\alpha _{1}+2\alpha
_{2}}-b_{3}c_{1}E_{\alpha _{1}+\alpha _{2}}-b_{4}c_{1}E_{\alpha _{1}}
\end{eqnarray}%
Therefore $[X,Y]\in \mathfrak{m}_{2\alpha }$. \qed

\begin{prop}
Let $G_2/U(2)$ be the flag manifold represented by the short root. Let $X\in%
\mathfrak{m}=T_oF$ be a nonzero vector. Then $X$ is equigeodesic if, and
only if, $X\in \mathfrak{m}_{\alpha_1}$ or $X\in\mathfrak{m}_{2\alpha_1}$.
\end{prop}

\noindent \textbf{Proof:}\quad Writing $X=X_{\alpha _{1}}+X_{2\alpha _{1}}$
with $X_{\alpha _{1}}\in \mathfrak{m}_{\alpha _{1}}$ and $X_{2\alpha
_{1}}\in \mathfrak{m}_{2\alpha _{1}}$, we have 
\begin{eqnarray}
\lbrack X,\Lambda X] &=&[X_{\alpha _{1}}+X_{2\alpha _{1}},\lambda
_{1}X_{\alpha _{1}}+\lambda _{2}X_{2\alpha _{1}}]  \notag \\
&=&\lambda _{1}[X_{\alpha _{1}},X_{\alpha _{1}}]+\lambda _{2}[X_{\alpha
_{1}},X_{2\alpha _{1}}]+\lambda _{1}[X_{2\alpha _{1}},X_{\alpha
_{1}}]+\lambda _{2}[X_{2\alpha _{1}},X_{2\alpha _{1}}]  \notag \\
&=&(\lambda _{2}-\lambda _{1})[X_{\alpha _{1}},X_{2\alpha _{1}}].
\end{eqnarray}%
If $X$ is equigeodesic then $(\lambda _{2}-\lambda _{1})[X_{\alpha
_{1}},X_{2\alpha _{1}}]=0$ for any $\lambda _{1}>0,\lambda _{2}>0$ and
therefore $[X_{\alpha _{1}},X_{2\alpha _{1}}]=0$. According the previous
lemma, $[X_{\alpha _{1}},X_{2\alpha _{1}}]=0$ if, and only if $X_{\alpha
_{1}}=0$ or $X_{2\alpha _{1}}=0$.

On other hand, if $X\in \mathfrak{m}_{\alpha_1}$ then $\Lambda X=\lambda_1 X$
for any invariant metric $\Lambda$ and $[X,\Lambda X]$=0 for any $\Lambda$
and $X$ is equigeodesic vector. Analogously for $X\in \mathfrak{m}%
_{2\alpha_1}$. \qed
\newline

\emph{b) Consider now $G_2(\alpha_1)=G_2/U(2)$ represented by the long root.}
In this case we have $R_K=\{ \alpha_1\}$ and $R_M=\{ \alpha_2 \}$. The
isotropic representation have three irreducible submodules and the tangent
space at $o$ slipts as $\mathfrak{m}=\mathfrak{m}_{\alpha_2}\oplus \mathfrak{%
m}_{2\alpha_2} \oplus \mathfrak{m}_{3\alpha_2}$, where $\mathfrak{m}%
_{\alpha_2}= \mathfrak{u}_{\alpha_2}\oplus \mathfrak{u}_{\alpha_1+\alpha_2}$%
, $\mathfrak{m}_{2\alpha_2}= \mathfrak{u}_{\alpha_1+2\alpha_2}$ and $%
\mathfrak{m}_{3\alpha_2}=\mathfrak{u}_{\alpha_1+3\alpha_2} \oplus \mathfrak{u%
}_{2\alpha_1+3\alpha_2}$. As in Lemma \ref{inclusao dois somandos} we prove

\begin{lema}
The following inclusions hold:

\begin{enumerate}
\item $[\mathfrak{m}_{\alpha_2},\mathfrak{m}_{2\alpha_2}]\subset \mathfrak{m}%
_{\alpha_2}\oplus \mathfrak{m}_{3\alpha_2}$;

\item $[\mathfrak{m}_{\alpha_2},\mathfrak{m}_{3\alpha_2}] \subset \mathfrak{m%
}_{2\alpha_2} $;

\item $[\mathfrak{m}_{2\alpha_2}, \mathfrak{m}_{3\alpha_2}]\subset \mathfrak{%
m}_{\alpha_2} $.
\end{enumerate}
\end{lema}

\noindent \textbf{Proof:}\quad Let $E_{\alpha}$ be a Weyl's basis of $%
\mathfrak{g}_2$. Writting 
\begin{equation}  \label{baseweyl}
\begin{array}{ccl}
X & = & a_1E_{\alpha_2}+a_2E_{\alpha_1+\alpha_2}+b_1E_{-\alpha_2}+b_2E_{-(%
\alpha_1+\alpha_2)} \in \mathfrak{m}_{\alpha_2}, \\ 
Y & = & c_1E_{\alpha_1+2\alpha_2} + c_2E_{-(\alpha_1+2\alpha_2)} \in 
\mathfrak{m}_{2\alpha_2}, \\ 
Z & = & d_1E_{\alpha_1+3\alpha_2}+d_2E_{2\alpha_1+3\alpha_2}+
e_1E_{-(\alpha_1+3\alpha_2)}+e_2E_{-(2\alpha_1+3\alpha_2)}\in \mathfrak{m}%
_{3\alpha_2},%
\end{array}%
\end{equation}

with $b_i=-\overline{a_i},c_2=-\overline{c_1}$ and $e_i=-\overline{d_i}$ we
have 
\begin{eqnarray*}
1)\, [X,Y]&=& a_1c_1E_{\alpha_1+3\alpha_2}+a_1c_2E_{-(\alpha_1+\alpha_2)}
+a_2c_1E_{2\alpha_1+3\alpha_2} + a_2c_2E_{-\alpha_2} \\
&-& b_1c_1E_{\alpha_1+\alpha_2}-b_1c_2E_{-(\alpha_1+3\alpha_2)} -
b_2c_1E_{\alpha_2} - b_2c_2E_{-(2\alpha_1+3\alpha_2)}\subset \mathfrak{m}%
_{\alpha_2}\oplus\mathfrak{m}_{3\alpha_2} ; \\
\\
2)\, [X,Z]&=&
a_1e_1E_{-(\alpha_1+2\alpha_2)}+a_2e_2E_{-(\alpha_1+2\alpha_2)} -
b_1d_1E_{\alpha_1+2\alpha_2}- b_2d_2E_{\alpha_1+2\alpha_2} \\
&=&
(a_1e_1+a_2e_2)E_{-(\alpha_1+2\alpha_2)}+(-b_1d_1-b_2d_2)E_{\alpha_1+2%
\alpha_2} \subset \mathfrak{m}_{2\alpha_2}; \\
\\
3)\, [Y,Z]&=&
c_1e_1E_{-\alpha_2}+c_1e_2E_{-(\alpha_1+\alpha_2)}-c_2d_1E_{%
\alpha_2}-c_2d_2E_{\alpha_1+\alpha_2} \subset \mathfrak{m}_{\alpha_2}.
\end{eqnarray*}
\qed
\newline

\begin{teo}
\label{g2alfa2} Consider the flag manifold $G_2(\alpha_2)$ and $V\in 
\mathfrak{m}$. Write $V=X+Y+Z$, where $X,Y,Z$ are as in (\ref{baseweyl}).
Then the equation $[V,\Lambda V]_{\mathfrak{m}}=0$, for any invariant metric 
$\Lambda$, is equivalent to the following system of algebraic equations 
\begin{equation}  \label{sistemaG2}
\left\{ 
\begin{array}{ccc}
b_2c_1 & = & 0 \\ 
c_2d_1 & = & 0 \\ 
b_1c_1 & = & 0 \\ 
c_2d_2 & = & 0 \\ 
b_1d_1+b_2d_2 & = & 0 \\ 
a_1c_1 & = & 0 \\ 
a_2c_1 & = & 0%
\end{array}
\right. .
\end{equation}

Therefore $V$ is an equigeodesic vector if, and only if, one of the
following holds: \newline
a) $V\in \mathfrak{m}_{\alpha _{2}}$;\newline
b) $V\in \mathfrak{m}_{2\alpha _{2}}$;\newline
c) $V\in \mathfrak{m}_{3\alpha _{2}}$;\newline
d) $V\in \mathfrak{u}_{\alpha _{1}+\alpha _{2}}\oplus \mathfrak{u}_{\alpha
_{1}+3\alpha _{2}}$;\newline
e) $c_{1}=0,c_{2}=0,d_{1}=-\frac{b2\ast d2}{b1}$.

\end{teo}

\noindent \textbf{Proof:}\quad Let $V=X+Y+Z \in \mathfrak{m}$, where $X,Y,Z$
are as in equation $(\ref{baseweyl})$. We have: 
\begin{eqnarray*}
[V,\Lambda V]_{\mathfrak{m}}&=& [X+Y+Z,\lambda_1X +\lambda_2Y+\lambda_3Z] \\
&=&
(\lambda_2-\lambda_1)[X,Y]+(\lambda_3-\lambda_1)[X,Z]+(\lambda_3-%
\lambda_2)[Y,Z] \\
&=& (\lambda_2-\lambda_1) \{
a_1c_1E_{\alpha_1+3\alpha_2}+a_1c_2E_{-(\alpha_1+\alpha_2)}
+a_2c_1E_{2\alpha_1+3\alpha_2} + a_2c_2E_{-\alpha_2} \\
&& -b_1c_1E_{\alpha_1+\alpha_2}-b_1c_2E_{-(\alpha_1+3\alpha_2)} -
b_2c_1E_{\alpha_2} - b_2c_2E_{-(2\alpha_1+3\alpha_2)} \} \\
&&
+(\lambda_3-\lambda_1)\{(a_1e_1+a_2e_2)E_{-(\alpha_1+2%
\alpha_2)}+(-b_1d_1-b_2d_2)E_{\alpha_1+2\alpha_2} \} \\
&& +(\lambda_3-\lambda_2)\{
c_1e_1E_{-\alpha_2}+c_1e_2E_{-(\alpha_1+\alpha_2)}-c_2d_1E_{%
\alpha_2}-c_2d_2E_{\alpha_1+\alpha_2} \}.
\end{eqnarray*}

Then $V$ is an equigeodesic vector if, and only if, the coefficients of $V$
satisfy the system of equations (\ref{sistemaG2}). 

The solutions of (\ref{sistemaG2}) are:\newline
a) $c_1 = 0, c_2 = 0, d_1 = -\frac{b2*d2}{b1}$;\newline
b) $c_1 = 0, d_1 = 0, d_2 = 0$; \newline
c) $b_1 = 0, c_1 = 0, c_2 = 0, d_2 = 0$; \newline
d) $a_1 = 0, a_2 = 0, b_1 = 0, b_2 = 0, c_2 = 0$; \newline
e) $a_1 = 0, a_2 = 0, b_1 = 0, b_2 = 0, d_1 = 0, d_2 = 0$; \newline
f) $b_1 = 0, b_2 = 0, c_1 = 0, c_2 = 0$.\newline

Therefore the solutions of the algebraic system (\ref{sistemaG2}) determine
all vector spaces that appear in the Theorem \ref{g2alfa2}.

\qed

\section{Rank three T-Roots and Einstein Metrics}

In the study of the geometry of flag manifolds a class of metrics play a key
role namely the Einstein metrics. 
We recall that an invariant metric $\Lambda $ on $\mathbb{F}$ is Einstein if
it is proportional to the Ricci tensor, that is, it satisfies $Ric_{\Lambda
}=c\Lambda $.


It is well known that the problem of finding invariant
Einstein metrics on flag manifolds reduces to solve a delicate algebraic system, \cite{wangziller}. A well known solution for this
system is the Kähler-Einstein one. Indeed, for each invariant complex
structure on $\mathbb{F}$ there exist a unique invariant Kähler-Einstein
metric, see (\cite{Besse}, 8.95).

There are several examples of invariant Einstein non-Kähler metrics. All of
them have repetition in the parameters of the metric, see \cite%
{Evan2} or \cite{neiton}. In this section we connect this repetition on the
parameters of an invariant Einstein metric with the dimension of the
isotropic summands $\mathfrak{m}_{i}$.

We derive it using the expression of
the scalar curvature correspondent to an invariant metric and our result
concerning zero sum triple of T-roots. We prove that there are no isolated
T-root, that is, every T-root belongs to a zero sum triple of T-roots.

Let $\left\langle \cdot ,\cdot \right\rangle $ be the Cartan-Killing form on 
$\mathfrak{g}$ and 
\begin{equation*}
\mathfrak{m}=\mathfrak{m}_{1}\oplus \ldots \oplus \mathfrak{m}_{s}
\end{equation*}
a decomposition into irreducible (non-equivalent) $\mbox{ad}(\mathfrak{k})-$%
submodules. Let $\{X_{\alpha }\}$ be a orthonormal basis (with respect to $%
-\left\langle \cdot ,\cdot \right\rangle $) adapted to the decomposition of $%
\mathfrak{m}$: $X_{\alpha }\in \mathfrak{m}_{i}$ and $X_{\beta }\in 
\mathfrak{m}_{j}$ with $i<j$ \ then $\alpha <\beta $. Following \cite%
{wangziller} we set $A_{\alpha \beta }^{\gamma }:=-\left\langle [X_{\alpha
},X_{\beta }],X_{\gamma }\right\rangle $, thus $[X_{\alpha },X_{\beta
}]=\sum_{\gamma }{A_{\alpha \beta }^{\gamma }X_{\gamma }}$. Consider 
\begin{equation}
C_{ij}^{k}:=\sum {(A_{\alpha \beta }^{\gamma })^{2}}  \label{coef cijk}
\end{equation}%
where the sum is taken over all indexes $\alpha ,\beta ,\gamma $ with $%
X_{\alpha }\in \mathfrak{m}_{i},X_{\beta }\in \mathfrak{m}_{j},X_{\gamma
}\in \mathfrak{m}_{k}$. Hence $C_{ij}^{k}$ is nonnegative, symmetric in all
of the three entries, and is independent of the orthonormal basis chosen for 
$\mathfrak{m}_{i},\mathfrak{m}_{j}$ and $\mathfrak{m}_{k}$ (but it depends
on the choice of the decomposition of $\mathfrak{m}$).

Let $\Lambda $ be an invariant metric on $\mathbb{F}$ and $S\left( \Lambda
\right) $ the correspondent scalar curvature. According to \cite{wangziller} 
\begin{equation}
S\left( \Lambda \right) =\frac{1}{2}\sum_{i}\frac{D_{i}}{\lambda _{i}}-\frac{%
1}{4}\sum_{i,j,k}C_{ij}^{k}\frac{\lambda _{k}}{\lambda _{i}\lambda _{j}}
\label{curvatura
escalar}
\end{equation}%
where $D_{i}=\dim _{\mathbb{R}}(\mathfrak{m}_{i})$ and $\lambda _{i}$
denotes the parameter of the invariant metric $\Lambda $ with $i=1,\ldots ,s$%
. We consider now the set of the invariant metrics with unitary volume: 
\begin{equation*}
\mathcal{M}=\{(\lambda _{1},\dots ,\lambda _{s})\in \mathbb{R}^{s}:\lambda
_{1}^{D_{1}}\cdots \lambda _{s}^{D_{s}}=1;\,\lambda _{1},\dots ,\lambda
_{s}>0\}.
\end{equation*}

The next result is in fact true for any compact, connected homogeneous
space. It shows an alternative manner of computing the Einstein
equations.
\begin{teo}\label{zw}(\cite{Besse}) Let
$\mathbb{F}$ be a flag manifold. Then the critical points of the
restriction map $S|_{\mathcal{M}}$ are precisely the invariant
Einstein metrics on $\mathbb{F}$.
\end{teo}

\begin{lema}(\cite{arva3})\label{lema soma}
Let $\xi,\eta,\zeta$ be T-roots such that $\xi+\eta+\zeta=0$. Then
there exist roots $\alpha,\beta,\gamma\in R$ with $k(\alpha)=\xi$,
$k(\beta)=\eta$, $k(\gamma)=\zeta$, such that
$\alpha+\beta+\gamma=0$.
\end{lema}

The calculus of the coefficients $C_{ij}^{k}$ can be laborius.
However the next result shows exactly which of them are nonzero.

%
\begin{lema}
\label{lema cijk} Let $\delta_{i},\delta_{j},\delta_{k}$ T-roots associated
to the real $ad(\mathfrak{k})$-modules $\mathfrak{m}_{i},\mathfrak{m}_{j}$
and $\mathfrak{m}_{k}$, respectively . Then $C_{ij}^{k}\neq0$ if, and only
if, $\delta_{i}+\delta_{j}+\delta_{k}=0$.
\end{lema}

\noindent \textbf{Proof:}\quad Consider the vectors $E_{\alpha }$, $\alpha
\in R_{M}$, of the fixed Weyl basis of $\mathfrak{g}^{\mathbb{C}}$ , and $%
V_{\alpha }:=\mathbb{R}S_{\alpha }+\mathbb{R}\sqrt{-1}A_{\alpha }$, $\alpha
\in R_{M}^{+}$. Hence, the vectors $I_{\alpha }=S_{\alpha }/\sqrt{2}$ and $%
F_{\alpha }=\sqrt{-1}A_{\alpha }/\sqrt{2}$, $\alpha \in R_{M}^{+}$ form a
orthonormal basis of $V_{\alpha }$. Thus, each set $b_{i}=\left\{ I_{\alpha
},F_{\alpha }:k(\alpha )=\delta _{i},\alpha \in R_{M}^{+}\right\} $ is a
orthonormal basis of a $ad(\mathfrak{k})$-module $\mathfrak{m}_{i}=\mathfrak{%
m}_{\delta _{i}}$, with $\delta _{i}\in R_{T}^{+}$.

Let $b_{i},b_{j}$ and $b_{k}$ be a orthonormal basis of $\mathfrak{m}_{i}$, $%
\mathfrak{m}_{j}$ and $\mathfrak{m}_{k}$, respectively. We notice that $%
\left[ \mathfrak{g}_{\alpha }^{\mathbb{C}},\mathfrak{g}_{\beta }^{\mathbb{C}}%
\right] =\mathfrak{g}_{\alpha +\beta }^{\mathbb{C}}$ and $\left( \mathfrak{g}%
_{\alpha +\beta }^{\mathbb{C}},\mathfrak{g}_{\gamma }^{\mathbb{C}}\right) =0$%
, unless $\alpha +\beta +\gamma =0$. Then $-\left\langle \left[ e_{\alpha
},e_{\beta }\right] ,e_{\gamma }\right\rangle =0$ except when $\alpha +\beta
+(-\gamma )=0$, for any $e_{\alpha }\in b_{i}$, $e_{\beta }\in b_{j}$ and $%
e_{\gamma }\in b_{k}$, with $\alpha ,\beta ,\gamma \in R_{M}^{+}$.

If $C_{ij}^{k}\neq 0$ then exists $\alpha ,\beta ,\gamma \in R_{M}^{+}$
with $k(\alpha )=\delta _{i}$, $k(\beta )=\delta _{j}$, $k(\gamma )=\delta
_{k}$ such that $\alpha +\beta +(-\gamma )=0$. Hence $\delta _{i}+\delta
_{j}+(-\delta _{k})=k\left( \alpha \right) +k(\beta )+(-k\left( \gamma
\right) )=k(\alpha +\beta +(-\gamma ))=0$.

Conversely, if $\delta _{i},\delta _{j},\delta _{k}$ are nonzero $T$-roots
such that $\delta _{i}+\delta _{j}+\delta _{k}=0$, then there exist $\alpha
,\beta ,\gamma \in R_{M}$ with $k(\alpha )=\delta _{i}$, $k(\beta )=\delta
_{j}$, $k(\gamma )=\delta _{k}$ such that $\alpha +\beta +\gamma =0$, hence $%
C_{ij}^{k}\neq 0$. \qed

\begin{defi}
Let $\mathbb{F}$ be a flag manifold and $R_{T}$ the correspondent set of
T-roots We say that a T-root \emph{$\delta _{i}$ belongs to a T-zero sum
triple} if there are T-roots $\delta _{j},\delta _{k}\in R_{T}$ such that $%
\delta _{i}+\delta _{j}+\delta _{k}=0$. In this case we denote by $T(\delta
_{i})$ the number of T-zero sum triple contained the T-root $\delta _{i}$.
Of course $T(\delta _{i})=T(-\delta _{i})$, for every $\delta _{i}\in R_{T}$.
\end{defi}

We now connect the repetition on the parameters of an invariant Einstein
metric with the dimension of its associated isotropic summands.

\begin{prop}
\label{teo da dim} Let\ $\delta _{i}$, $\delta _{j}$, $\delta _{k}\in $ $%
R_{T}$ with $\delta _{i}$, $\delta _{j}$, $\delta _{k}=-(\delta _{i}+\delta
_{j})$ such that $T(\delta _{i})=T(\delta _{j})=1$. If there exists $i$ and $%
j$ such that $\mathbb{F}$ admits an invariant Einstein metric $\Lambda $
satisfying $\lambda _{i}=\lambda _{j}$ then $\dim \mathfrak{m}_{i}=\dim 
\mathfrak{m}_{j}$.
\end{prop}

\noindent \textbf{Proof:}\quad According Theorem \ref{zw} an invariant
metric $\Lambda $ is Einstein if, and only if, $\Lambda $ is solution of the 
$s+1$ equations 
\begin{equation}
\frac{\partial S}{\partial \lambda _{l}}=\xi D_{l}\lambda _{1}^{D_{1}}\cdots
\lambda _{l}^{D_{l}-1}\cdots \lambda _{s}^{D_{s}},\text{ }1\leq l\leq s
\label{eq Einstein por S}
\end{equation}%
\begin{equation}
\lambda _{1}^{D_{1}}\cdots \lambda _{s}^{D_{s}}=1  \label{vol}
\end{equation}%
where $D_{l}=\dim \mathfrak{m}_{l}$ and $\xi $ denotes the Lagrange
multiplier.

In particular, $\Lambda $ must satisfy the two equations of (\ref{eq
Einstein por S}) for $l=i$ and $l=j$. By hypothesis $T(\delta _{i})=T(\delta
_{j})=1$, then using Lemma \ref{coef cijk} and formula (\ref{curvatura
escalar}) we conclude that the equation ($\ref{eq Einstein por S}$), for $%
l=i $, reduces to 
\begin{equation}
-\frac{D_{i}}{2\lambda _{i}^{2}}-\frac{1}{4}C_{ij}^{k}\left( \frac{1}{%
\lambda _{j}\lambda _{k}}-\frac{\lambda _{k}}{\lambda _{i}^{2}\lambda _{j}}-%
\frac{\lambda _{j}}{\lambda _{i}^{2}\lambda _{k}}\right) =\xi D_{i}\lambda
_{1}^{D_{1}}\cdots \lambda _{i}^{D_{i}-1}\cdots \lambda _{s}^{D_{s}}.
\end{equation}%
Multiplying this equality by $\lambda _{i}/D_{i}$ and using (\ref{vol}) we
obtain 
\begin{equation}
-\frac{1}{2\lambda _{i}}-\frac{1}{4D_{i}}C_{ij}^{k}\left( \frac{\lambda _{i}%
}{\lambda _{j}\lambda _{k}}-\frac{\lambda _{k}}{\lambda _{i}\lambda _{j}}-%
\frac{\lambda _{j}}{\lambda _{i}\lambda _{k}}\right) =\xi .  \label{lambdai}
\end{equation}

Analogously, if $l=j$ we obtain 
\begin{equation}
-\frac{1}{2\lambda_{j}}-\frac{1}{4D_{j}}C_{ij}^{k}\left(\frac{\lambda_{j}}{%
\lambda_{i}\lambda_{k}} -\frac{\lambda_{k}}{\lambda_{i}\lambda_{j}}-\frac{%
\lambda_{i}}{\lambda_{j}\lambda_{k}}\right)=\xi.  \label{lambdaj}
\end{equation}

Therefore if $\Lambda $ is an invariant Einstein metric satisfying $\lambda
_{i}=\lambda _{j}$, according to equations (\ref{lambdai}) and (\ref{lambdaj}%
), we obtain $D_{i}=D_{j}$. \qed\newline

The next result shows that the set of $T$-roots enjoys an interesting
property, despite not being a root system.


\begin{teo}
\label{soma tripla} Let $\mathbb{F}$ be a flag manifold and $R_{T}$ a
correspondent set of T-roots. If $R_{T}$ contain more than one positive
T-root, then every T-root belongs to some T-zero sum triple.
\end{teo}

\noindent \textbf{Proof:}\quad Consider the invariant Kähler-Einstein metric 
$\Lambda _{K}$ on $\mathbb{F}$ associated to the canonical invariant complex
structure on $\mathbb{F}$. $\Lambda _{K}$ satisfies the algebraic system (%
\ref{eq Einstein por S}) and (\ref{vol}).

Suppose the existence of a T-root $\delta _{k_{0}}$ which do not belong to
any T-zero sum triple. According to Lemma \ref{lema cijk}, we see that $%
C_{ij}^{k_{0}}=0$ for every $i,j=1,\dots ,s$. Thus, the equation (\ref{eq
Einstein por S}) for $l=k_{0}$ reduces to 
\begin{equation*}
-\frac{D_{k_{0}}}{2\lambda _{k_{0}}^{2}}=\xi D_{k_{0}}\lambda
_{1}^{D_{1}}\cdots \lambda _{k_{0}}^{D_{k_{0}}-1}\cdots \lambda _{s}^{D_{s}}
\end{equation*}%
then 
\begin{equation}
\xi =-\frac{1}{2\lambda _{k_{0}}}.
\end{equation}

By assumption there exist a positive T-root $\delta_{i}\neq \delta_{k_{0}}$.
Then the equation (\ref{eq Einstein por S}) for $l=i$ becomes 
\begin{equation*}
\frac{\partial S}{\partial \lambda_{i}}=-\frac{1}{2\lambda_{k_{0}}}%
D_{i}\lambda_{1}^{D_{1}}\cdots\lambda_{i}^{D_{i}-1}\cdots\lambda_{s}^{D_{s}}.
\end{equation*}
Then, 
\begin{equation}
\frac{\partial S}{\partial\lambda_{i}}=-\frac{D_{i}}{2\lambda_{i}%
\lambda_{k_{0}}}.  \label{derivada}
\end{equation}

On the other hand, 
\begin{equation}
\frac{\partial S}{\partial \lambda _{i}}=-\frac{D_{i}}{2\lambda _{i}^{2}}+f
\label{deriv
2}
\end{equation}%
, $f$ \ being a function depending on each $\lambda _{j}$ whose T-roots
satisfies $\delta _{i}\pm \delta _{j}\in R_{T}$ ( according to \ref{lema
cijk} ). Therefore using (\ref{derivada}) and (\ref{deriv 2}) we see that

\begin{equation*}
f=\frac{D_{i}}{2\lambda_{i}}\left(\frac{1}{\lambda_{i}}-\frac{1}{%
\lambda_{k_{0}}}\right).
\end{equation*}
This equality contradicts the fact that $\delta_{k_{0}}$ do not belongs to
any T-zero sum triple. This concludes the proof. \qed

This result allow us to characterize some sets of T-roots.

\begin{corol}
\label{duas t-raizes positivas} Let $\mathbb{F}$ be a flag manifold such
that $\mathfrak{m}=\mathfrak{m}_{1}\oplus\mathfrak{m}_{2}$, then $%
R_{T}=\{\pm \zeta,\pm 2\zeta\}$ where $\zeta\in\mathfrak{t}^{*}\setminus{0}$%
. Also, if $\dim\mathfrak{m}_{1}\neq \dim\mathfrak{m}_{2}$ then any
invariant Einstein metric on $\mathbb{F}$ satisfies $\lambda_{1}\neq%
\lambda_{2}$.
\end{corol}

\noindent \textbf{Proof:}\quad As $\mathfrak{m}=\mathfrak{m}_{1}\oplus 
\mathfrak{m}_{2}$ we may write the set of T-roots in the form $R_{T}=\left\{
\delta ,-\delta ,\zeta ,-\zeta \right\} $ with $\delta ,\zeta \in k\left(
R^{+}\right) $. By Theorem \ref{soma tripla}, $\delta $ belongs to some
T-zero sum triple. But, T-roots are nonzero linear functionals in $\mathfrak{t%
}^{\ast }$, then the possibilities for the T-zero sum triple containing $%
\delta $ are $\delta +\zeta +\zeta =0$, $\delta +\delta +\zeta =0$ and $%
\delta -\zeta -\zeta =0$. According to Lemma \ref{lema soma}, we see that
the first and the second possibilities cannot happens because $\delta ,\zeta
\in k\left( R^{+}\right) $, then $\delta =2\zeta $.

For the remaining possibility we see that each T-root belongs to exactly one
T-zero sum triple. Now Proposition \ref{teo da dim} follows from the
Corollary. \qed

In the case of three isotropic summands (irreducibles and inequivalent) some
T-roots may belong to more than one T-zero sum triple, but even in this case
we can see that there exist few possibilities for the set of T-roots.

\begin{corol}
\label{corol3somandos} Let $\mathbb{F}$ be a flag manifold such that $%
\mathfrak{m}=\mathfrak{m}_{1}\oplus\mathfrak{m}_{2}\oplus\mathfrak{m}_{3}$,
then $R_{T}=\{\pm \delta,\pm \zeta,\pm(\delta+\zeta)\}$ or $R_{T}=\{\pm
\delta,\pm 2\delta,\pm 4\delta\}$, where $\delta,\zeta\in\mathfrak{t}%
^{*}\setminus \{0\}$.
\end{corol}

\noindent \textbf{Proof:}\quad By assumption, the set of T-roots is given by 
$R_{T}=\left\{ \pm \alpha ,\pm \beta ,\pm \gamma \right\} $ where $\alpha
,\beta ,\gamma \in k(R^{+})$ with $\alpha ,\beta ,\gamma \in \mathfrak{t}%
^{\ast }\setminus \{0\}$.

But do not exist T-zero sum triple with only positive T-roots or containing
two opposite sign T-roots. Thus, the possibilities for the T-zero sum
triples containing $\alpha $ are $\left( \alpha ,\beta ,-\gamma \right) $, $%
\left( \alpha ,-\beta ,\gamma \right) $, $\left( \alpha ,-\beta ,-\gamma
\right) $, $\left( \alpha ,\alpha ,-\beta \right) $, $\left( \alpha ,\alpha
,-\gamma \right) $, $\left( \alpha ,-\beta ,-\beta \right) $, $\left( \alpha
,-\gamma ,-\gamma \right) $.

For any of the first three choices we conclude that the set of T-roots has
the form $R_{T}=\{\pm \delta ,\pm \zeta ,\pm (\delta +\zeta )\}$. For the
last four choices we obtain that the set of T-roots has the form $%
R_{T}=\{\pm \delta ,\pm 2\delta ,\pm 3\delta \}$ or $R_{T}=\{\pm \delta ,\pm
2\delta ,\pm 4\delta \}$. \qed

\begin{exemplo}
\textrm{\ Consider the flag manifold }$G_{2}(\alpha _{2})$.\textrm{It is a direct computation to verify that in this case $%
\mathfrak{m}=\mathfrak{m}_{1}\oplus \mathfrak{m}_{2}$ with $\dim \mathfrak{m}%
_{1}=8$ and $\dim \mathfrak{m}_{2}=2$. Hence, any invariant Einstein metric
has exactly two parameters $\lambda _{1}$ and $\lambda _{2}$ with $\lambda
_{1}\neq \lambda _{2}$, according to Corollary \ref{corol3somandos}.} \qed
\end{exemplo}

\begin{exemplo}
\textrm{{\ Consider the flag manifold $\mathbb{F}=$}}$G_{2}(\alpha _{1})$%
\textrm{{. In this case the subalgebra $\mathfrak{t}$ has the form 
\begin{equation*}
diag\{\delta ,\delta ,-2\delta \}\in \mathfrak{su}(3)
\end{equation*}%
and a choice of positive T-roots is the same as a choice of a set of linear
functionals of the form $R_{T}^{+}=\{\delta ,2\delta ,3\delta \}$. Then $%
\mathfrak{m}=\mathfrak{m}_{1}\oplus \mathfrak{m}_{2}\oplus \mathfrak{m}_{3}$
where $\mathfrak{m}_{1},\mathfrak{m}_{2}$ and $\mathfrak{m}_{3}$ correspond
to the functionals $\delta ,2\delta $ and $3\delta $, respectively. }}

\textrm{It is easy to check that the T-roots $\delta $ and $2\delta $ belong
to exactly one T-zero sum triple and $\dim \mathfrak{m}_{1}=4$, $\dim 
\mathfrak{m}_{2}=2$. Then any invariant Einstein metric on $\mathbb{F}$ has
exactly three parameters $\lambda _{1},\lambda _{2}$ and $\lambda _{3}$
satisfying $\lambda _{1}\neq \lambda _{2}$. As $\dim \mathfrak{m}_{3}=4$ and 
$T(3\delta )=1$. In a similar way, we prove that any invariant Einstein
metric must satisfy $\lambda _{1}\neq \lambda _{3}$. This gives an
alternative proof of the well known fact that the normal metric is not
Einstein on $G_{2}/U(2)$. \qed}
\end{exemplo}

We, just for completeness, conclude this section determining the Einstein
equations for the full flag manifold $G_{2}/T$ and derive explicitly, the
invariant Kähler-Einstein metric corresponding to each invariant complex
structure.

In the case of maximal flag manifolds, the irreducible and inequivalent
isotropic summands $\mathfrak{m}_{i}$ are determined by the positive roots $%
R^{+}$. So they are indexed by these roots and for the coefficients $%
C_{ij}^{k}$ we may write $C_{\alpha \beta }^{\gamma }$ with $\alpha ,\beta
,\gamma \in R^{+}$.

According to \cite{sakane} if $\mathbb{F}=G/T$ is a full flag manifold
then the Einstein equation of a invariant metric $g$ on $\mathbb{F}$ is
given by

\begin{equation}
\label{r alpha}
c=\frac{1}{2\lambda _{\alpha }}+\frac{1}{8}\sum_{\beta ,\gamma \in R^{+}}%
\frac{\lambda _{\alpha }}{\lambda _{\beta }\lambda _{\gamma }}C_{\beta
\gamma }^{\alpha }-\frac{1}{4}\sum_{\beta ,\gamma \in R^{+}}\frac{\lambda
_{\gamma }}{\lambda _{\alpha }\lambda _{\beta }}C_{\alpha \beta }^{\gamma }
\end{equation}%
where $g(\cdot ,\cdot )=\lambda _{\alpha }(\cdot ,\cdot )|_{\mathfrak{m}%
_{\alpha }},\text{ }\alpha \in R^{+}$, $c$ is the Einstein constant and $T$
is a maximal torus on $G$.

Now, using the fixed Weyl base of $\mathfrak{g}^{\mathbb{C}}$ we see that
the unique triple of positive roots in $\mathfrak{g}^{\mathbb{C}}$ such that 
$C_{\alpha \beta }^{\gamma }$ is nonzero are

\begin{equation}
C_{\alpha \beta }^{\alpha +\beta }=2\left( N_{\alpha ,\beta }\right) ^{2}%
\hspace{0.5cm}\text{ and }\hspace{0.5cm}C_{\alpha \beta }^{\alpha -\beta
}=2\left( N_{\alpha ,-\beta }\right) ^{2}  \label{colchetes g2}
\end{equation}%
where $N_{\alpha ,\beta }$ are the constants of structure of $\mathfrak{g}^{%
\mathbb{C}}$. Thus, we obtain the Einstein equation for $G_{2}/T$.

\begin{prop}
The Einstein equations for the full flag manifold $\mathbb{F}=G_{2}/T$ are
given by the following algebraic system 
\begin{eqnarray*}
c&=& \frac{1}{2\lambda_{\alpha_1}}+\frac{1}{16}\left( \frac{%
\lambda_{\alpha_1}}{\lambda_{\alpha_1+\alpha_2}\lambda_{\alpha_2}}+\frac{%
\lambda_{\alpha_1}}{\lambda_{2\alpha_1+3\alpha_2}\lambda_{\alpha_1+3\alpha_2}%
}\right)-\frac{1}{16}\left( \frac{\lambda_{\alpha_2}}{\lambda_{\alpha_1}%
\lambda_{\alpha_1+\alpha_2}}+\frac{\lambda_{\alpha_1+\alpha_2}}{%
\lambda_{\alpha_1}\lambda_{\alpha_2}}+\frac{\lambda_{\alpha_1+3\alpha_2}}{%
\lambda_{\alpha_1}\lambda_{2\alpha_1+3\alpha_2}}+\frac{\lambda_{2\alpha_1+3%
\alpha_2}}{\lambda_{\alpha_1}\lambda_{\alpha_1+3\alpha_2}}\right) \\
c&=& \frac{1}{2\lambda_{\alpha_2}}+\frac{1}{16}\left( \frac{%
\lambda_{\alpha_2}}{\lambda_{\alpha_1+\alpha_2}\lambda_{\alpha_1}}+\frac{%
\lambda_{\alpha_2}}{\lambda_{\alpha_1+3\alpha_2}\lambda_{\alpha_1+2\alpha_2}}%
\right) +\frac{1}{12}\frac{\lambda_{\alpha_2}}{\lambda_{\alpha_1+2\alpha_2}%
\lambda_{\alpha_1+\alpha_2}} -\frac{1}{12}\left( \frac{\lambda_{\alpha_1+%
\alpha_2}}{\lambda_{\alpha_2}\lambda_{\alpha_1+2\alpha_2}}+\frac{%
\lambda_{\alpha_1+2\alpha_2}}{\lambda_{\alpha_2}\lambda_{\alpha_1+\alpha_2}}%
\right) \\
&& -\frac{1}{16}\left(\frac{\lambda_{\alpha_1}}{\lambda_{\alpha_2}\lambda_{%
\alpha_1+\alpha_2}}+\frac{\lambda_{\alpha_1+\alpha_2}}{\lambda_{\alpha_2}%
\lambda_{\alpha_1}}+\frac{\lambda_{\alpha_1+2\alpha_2}}{\lambda_{\alpha_2}%
\lambda_{\alpha_1+3\alpha_2}}+\frac{\lambda_{\alpha_1+3\alpha_2}}{%
\lambda_{\alpha_2}\lambda_{\alpha_1+2\alpha_2}}\right) \\
c&=&\frac{1}{2\lambda_{\alpha_1+\alpha_2}}+\frac{1}{16}\left( \frac{%
\lambda_{\alpha_1+\alpha_2}}{\lambda_{\alpha_1}\lambda_{\alpha_2}}+\frac{%
\lambda_{\alpha_1+\alpha_2}}{\lambda_{2\alpha_1+3\alpha_2}\lambda_{%
\alpha_1+2\alpha_2}}\right) +\frac{1}{12}\frac{\lambda_{\alpha_1+\alpha_2}}{%
\lambda_{\alpha_1+2\alpha_2}\lambda_{\alpha_2}} -\frac{1}{12}\left( \frac{%
\lambda_{\alpha_2}}{\lambda_{\alpha_1+\alpha_2}\lambda_{\alpha_1+2\alpha_2}}+%
\frac{\lambda_{\alpha_1+2\alpha_2}}{\lambda_{\alpha_1+\alpha_2}\lambda_{%
\alpha_2}}\right) \\
&&-\frac{1}{16}\left( \frac{\lambda_{\alpha_2}}{\lambda_{\alpha_1+\alpha_2}%
\lambda_{\alpha_1}}+\frac{\lambda_{\alpha_1}}{\lambda_{\alpha_1+\alpha_2}%
\lambda_{\alpha_2}}+\frac{\lambda_{\alpha_1+2\alpha_2}}{\lambda_{\alpha_1+%
\alpha_2}\lambda_{2\alpha_1+3\alpha_2}}+\frac{\lambda_{2\alpha_1+3\alpha_2}}{%
\lambda_{\alpha_1+\alpha_2}\lambda_{\alpha_1+2\alpha_2}}\right) \\
c&=&\frac{1}{2\lambda_{\alpha_1+2\alpha_2}}+\frac{1}{16}\left( \frac{%
\lambda_{\alpha_1+2\alpha_2}}{\lambda_{\alpha_1+3\alpha_2}\lambda_{\alpha_2}}%
+\frac{\lambda_{\alpha_1+2\alpha_2}}{\lambda_{2\alpha_1+3\alpha_2}\lambda_{%
\alpha_1+\alpha_2}}\right) +\frac{1}{12}\frac{\lambda_{\alpha_1+2\alpha_2}}{%
\lambda_{\alpha_1+\alpha_2}\lambda_{\alpha_2}} -\frac{1}{12}\left( \frac{%
\lambda_{\alpha_2}}{\lambda_{\alpha_1+2\alpha_2}\lambda_{\alpha_1+\alpha_2}}+%
\frac{\lambda_{\alpha_1+\alpha_2}}{\lambda_{\alpha_1+2\alpha_2}\lambda_{%
\alpha_2}}\right) \\
&& -\frac{1 }{16}\left( \frac{\lambda_{\alpha_2}}{\lambda_{\alpha_1+2%
\alpha_2}\lambda_{\alpha_1+3\alpha_2}}+\frac{\lambda_{\alpha_1+3\alpha_2}}{%
\lambda_{\alpha_1+2\alpha_2}\lambda_{\alpha_2}}+\frac{\lambda_{\alpha_1+%
\alpha_2}}{\lambda_{\alpha_1+2\alpha_2}\lambda_{2\alpha_1+3\alpha_2}}+\frac{%
\lambda_{2\alpha_1+3\alpha_2}}{\lambda_{\alpha_1+2\alpha_2}\lambda_{%
\alpha_1+\alpha_2}}\right) \\
c&=&\frac{1}{2\lambda_{\alpha_1+3\alpha_2}} -\frac{1}{16}\left( \frac{%
\lambda_{\alpha_2}}{\lambda_{\alpha_1+3\alpha_2}\lambda_{\alpha_1+2\alpha_2}}%
+\frac{\lambda_{\alpha_1+2\alpha_2}}{\lambda_{\alpha_1+3\alpha_2}\lambda_{%
\alpha_2}}+\frac{\lambda_{\alpha_1}}{\lambda_{\alpha_1+3\alpha_2}\lambda_{2%
\alpha_1+3\alpha_2}}+\frac{\lambda_{2\alpha_1+3\alpha_2}}{%
\lambda_{\alpha_1+3\alpha_2}\lambda_{\alpha_1}}\right) \\
&& +\frac{1}{16}\left( \frac{\lambda_{\alpha_1+3\alpha_2}}{%
\lambda_{\alpha_1+2\alpha_2}\lambda_{\alpha_2}}+\frac{\lambda_{\alpha_1+3%
\alpha_2}}{\lambda_{2\alpha_1+3\alpha_2}\lambda_{\alpha_1}}\right) \\
c&=&\frac{1}{2\lambda_{2\alpha_1+3\alpha_2}} -\frac{1}{16}\left( \frac{%
\lambda_{\alpha_1+3\alpha_2}}{\lambda_{2\alpha_1+3\alpha_2}\lambda_{\alpha_1}%
}+\frac{\lambda_{\alpha_1}}{\lambda_{2\alpha_1+3\alpha_2}\lambda_{\alpha_1+3%
\alpha_2}}+\frac{\lambda_{\alpha_1+2\alpha_2}}{\lambda_{2\alpha_1+3\alpha_2}%
\lambda_{\alpha_1+\alpha_2}}+\frac{\lambda_{\alpha_1+\alpha_2}}{%
\lambda_{2\alpha_1+3\alpha_2}\lambda_{\alpha_1+2\alpha_2}}\right) \\
&& +\frac{1}{16}\left( \frac{\lambda_{2\alpha_1+3\alpha_2}}{%
\lambda_{\alpha_1}\lambda_{\alpha_1+3\alpha_2}}+\frac{\lambda_{2\alpha_1+3%
\alpha_2}}{\lambda_{\alpha_1+\alpha_2}\lambda_{\alpha_1+2\alpha_2}}\right).
\end{eqnarray*}
\end{prop}

\dem Using (\ref{colchetes g2}) we compute $C_{\alpha\beta}^{\gamma}$, with $\alpha,\beta,\gamma\in\R^{+}$. Now the result follows from equation (\ref{r alpha}).\qed
\\

%


Let $\mathbb{F}$ be a flag manifold with a invariant complex structure $J$
fixed. We recall that there exists a bijection between partial ordering in $%
R_{M}$ and complex structures on $\mathbb{F}$. This correspondence is given
by 
\begin{equation*}
JE_{\alpha }=\pm iE_{\alpha }\hspace{0.5cm}\alpha \in R_{M}^{+}.
\end{equation*}%
%


It is well known that for each invariant complex structure (or equivalently,
partial ordering in $R_{M}$) there exist a unique (up to homotheties)
invariant Kähler-Einstein metric, see (\cite{Besse},Chapter 8) or (\cite%
{sakane}). This metric is given by 
\begin{equation}
\Lambda _{J}=\{\lambda _{\alpha }=(\delta ,\alpha ):\delta =\frac{1}{2}%
\sum_{\beta \in R_{M}^{+}}\beta \}  \label{formulakahler}
\end{equation}%
where $(\cdot ,\cdot )$ is a inner product on $\mathfrak{h}^{\ast }$ induced
by the Cartan-Killing form of $\mathfrak{g}$.

According to \cite{borel} given a invariant complex structure $J$ on $%
\mathbb{F}$, we have a simple criterion satisfied for invariant Kähler
metrics on flags: a invariant metric $\Lambda $ is Kähler (with respect to $%
J $) if and only if 
\begin{equation}
\lambda _{\alpha +\beta }=\lambda _{\alpha }+\lambda _{\beta }\hspace{0.5cm}%
\alpha ,\beta \in R_{M}^{+}.  \label{mk}
\end{equation}%
%
%

But, for each Weyl chamber of the usual root system of $\mathfrak{g}_{2}$ we
have a choice of positive roots. As this root system has twelve roots and
they form successive angles of 30º, there are exactly twelve Weyl chambers,
then we have twelve possible invariant complex structure or six pairs of
conjugate structures.

On the other hand, if an invariant metric $\Lambda $ is Kähler with respect
to $J$, then $\Lambda $ is also Kähler with respect to $-J$. So it
sufficient to consider only the non conjugate invariant complex structure to
describe all the invariant Kähler-Einstein metrics. For each invariant
complex structure we have a choice of positive roots.

Now we describe explicitly the (unique) invariant Kahler-Einstein metric
associated to each complex structure on $G_{2}/T$.


\begin{prop}
The maximal flag manifold $G_2/T$ admits exactly (up to homotheties) six
invariant Kähler-Einstein metrics. These metrics and the correspond choice
of positive roots are given in the following table.

\begin{equation*}
\begin{tabular}{|c|c|c|}
\hline
$\Lambda$ & $R^{+}$ & $\Sigma$ \\ \hline
$(3,1,4,5,6,9)$ & $\alpha_{2},\alpha_{1}+3\alpha_{2},\alpha_{1}+2\alpha
_{2},2\alpha_{1}+3\alpha_{2},\alpha_{1}+\alpha_{2},\alpha_{1}$ & $\alpha
_{1},\alpha_{2}$ \\ \hline
$\left( {\ 6,5,1,4,9,3}\right) $ & $-(\alpha_{1}+\alpha_{2}),-\alpha_{1},%
\alpha
_{2},\alpha_{1}+3\alpha_{2},\alpha_{1}+2\alpha_{2},2\alpha_{1}+3\alpha_{2}$
& $-(\alpha_{1}+\alpha_{2}),2\alpha_{1}+3\alpha_{2}$ \\ \hline
$(3,4,1,5,9,6)$ & $-\alpha_{1},\alpha_{2},\alpha_{1}+3\alpha_{2},%
\alpha_{1}+2\alpha _{2},2\alpha_{1}+3\alpha_{2},\alpha_{1}+\alpha_{2}$ & $%
-\alpha_{1},\alpha _{1}+\alpha_{2}$ \\ \hline
$(6,1,5,4,3,9)$ & $\alpha_{1}+3\alpha_{2},\alpha_{1}+2\alpha_{2},2%
\alpha_{1}+3\alpha _{2},\alpha_{1}+\alpha_{2},\alpha_{1},-\alpha_{2}$ & $%
\alpha_{1}+3\alpha _{2},-\alpha_{2}$ \\ \hline
$(9,4,5,1,3,6)$ & $\alpha_{1}+2\alpha_{2},2\alpha_{1}+3\alpha_{2},%
\alpha_{1}+\alpha_{2},\alpha_{1},-\alpha_{2},-(\alpha_{1}+3\alpha_{2})$ & $%
\alpha_{1}+2\alpha _{2},-(\alpha_{1}+3\alpha_{2})$ \\ \hline
$(9,5,4,1,6,3)$ & $2\alpha_{1}+3\alpha_{2},\alpha_{1}+\alpha_{2},%
\alpha_{1},-\alpha _{2},-(\alpha_{1}+3\alpha_{2}),-(\alpha_{1}+2\alpha_{2})$
& $2\alpha _{1}+3\alpha_{2},-(\alpha_{1}+2\alpha_{2})$ \\ \hline
\end{tabular}%
\end{equation*}
where in the first column $\Lambda=(\lambda_{\alpha_1},\lambda_{\alpha_2},%
\lambda_{\alpha_1+\alpha_2},\lambda_{\alpha_1+2\alpha_2},\lambda_{\alpha_1+3%
\alpha_2},\lambda_{2\alpha_1+3\alpha_2}).$
\end{prop}

\noindent \textbf{Proof:}\quad The proof is obtained using the formula (\ref%
{formulakahler}). We perform the canonical choice for positives roots
(corresponding to the first row of the above table). The proof for the other
rows in the table are done in similar way. The canonical choice for the
roots is $R^+=\{\alpha_1,\alpha_2,\alpha_1+\alpha_2,\alpha_1+2\alpha_2,%
\alpha_1+3\alpha_2,2\alpha_1+3\alpha_2\}$.



Consider the fundamental weights $\Lambda_1$, $\Lambda_2$ related to simple
roots $\alpha_1,\alpha_2$ and defined by $\frac{2(\Lambda_i,\alpha_j)}{%
(\alpha_j,\alpha_j)}=\delta_{ij}$, $i,j=1,2$. Using the Cartan matrix of $%
\mathfrak{g}_2$ we can write the simple roots in terms of the fundamental
weights, see \cite{H}. The Cartan matrix of $\mathfrak{g}_2$ is given by 
\begin{equation*}
\left( 
\begin{matrix}
2 & -1 \\ 
-3 & 2%
\end{matrix}
\right),
\end{equation*}
and therefore $\alpha_1=2\Lambda_1-3\Lambda_2$ and $\alpha_2=-\Lambda_1+2%
\Lambda_2$. Then,

\begin{eqnarray*}
2\delta &=&
\alpha_1+\alpha_2+\alpha_1+\alpha_2+\alpha_1+2\alpha_1+\alpha_1+3\alpha_2+2%
\alpha_1+3\alpha_2 \\
&=& 6\alpha_1 +10\alpha_2 \\
&=& 6(2\Lambda_1-3\Lambda_2)+10(-\Lambda_1+2\Lambda_2) \\
&=& 2\Lambda_1+2\Lambda_2,
\end{eqnarray*}
and $\delta=\Lambda_1+\Lambda_2$. But $\frac{%
2(\Lambda_i,\alpha_j)}{(\alpha_j,\alpha_j)}=\delta_{ij}$ and setting $%
(\alpha_1,\alpha_1)=(\alpha_1+3\alpha_2,\alpha_1+3\alpha_2)=(2\alpha_1+3%
\alpha_2,2\alpha_1+3\alpha_2)=3$ (for long roots) and $(\alpha_2,\alpha_2)=(%
\alpha_1+\alpha_2,\alpha_1+\alpha_2)=(\alpha_1+2\alpha_2,\alpha_1+2%
\alpha_2)=1$ (for short roots) we have 
\begin{eqnarray*}
\lambda_{\alpha_1}&=& (\Lambda_1+\Lambda_2,\alpha_1)=\frac{1}{2}%
(\alpha_1,\alpha_1)=\frac{3}{2} \\
\lambda_{\alpha_2}&=& (\Lambda_1+\Lambda_2,\alpha_2)=(\Lambda_2,\alpha_2)=%
\frac{1}{2}(\alpha_2,\alpha_2)=\frac{1}{2} \\
\lambda_{\alpha_1+\alpha_2} &=&
(\Lambda_1+\Lambda_2,\alpha_1+\alpha_2)=(\Lambda_1,\alpha_1)+(\Lambda_2,%
\alpha_2)=2 \\
\lambda_{\alpha_1+2\alpha_2}&=&(\Lambda_1+\Lambda_2,\alpha_1+2\alpha_2)=(%
\Lambda_1,\alpha_1)+2(\Lambda_2,\alpha_2)=\frac{5}{2} \\
\lambda_{\alpha_1+3\alpha_2}&=&
(\Lambda_1+\Lambda_2,\alpha_1+3\alpha_2)=(\Lambda_1,\alpha_1)+3(\Lambda_2,%
\alpha_2)=3 \\
\lambda_{2\alpha_1+3\alpha_2}&=&
(\Lambda_1+\Lambda_2,2\alpha_1+3\alpha_2)=2(\Lambda_1,\alpha_1)+3(\Lambda_2,%
\alpha_2)= \frac{9}{2}.
\end{eqnarray*}
Thus an invariant K\"ahler-Einstein metric on $G_2/T$ is given (up to scale)
by $\Lambda=(3/2,1/2,2,5/2,3,9/2)$, and after normalization, we obtain the
metric in the table. \qed


\section*{Acknowledgments}

We would like to thank Luiz A.B. San Martin and Nir Cohen for helpful
discussions and sugestions. We also want to thank Ricardo Miranda Martins
for many useful comments. This research was partially supported by CNPq
grant 142271/2005-5 (da Silva), CAPES/CNPq grant 140431/2009-8 (Grama) and
FAPESP grant 07/06896-5 (Negreiros).


\begin{thebibliography}{99}
\addcontentsline{toc}{section}{References}

\bibitem{arva3} D. Alekseevsky and A. Arvanitoyeorgos; {Riemannian flag
manifolds with homogeneous geodesics,} \emph{Trans. Amer. Math. Soc.} 
\textbf{359} (2007), 1117--1126.

\bibitem{ale1} D. Alekseevsky; {Isotropy representation of flag manifolds, } 
\emph{Rend. Circ. Mat. Palermo (2) Suppl.} \textbf{54} (1998), 13--24.

\bibitem{Alek e Perel} D. V. Alekseevsky and A. M. Perelomov, \textit{%
Invariant Kähler-Einstein metrics on compact homogeneous spaces}, Funct.
Anal. Appl. 20 (1986),171-182.



\bibitem{Besse} A.L. Besse, \textit{Einstein Manifolds}, Springer-Verlag
1987.

\bibitem{borel} A.Borel, \emph{K\"ahlerian coset spaces of semisimple Lie
groups, } Proc. Nat. Acad. Sci. U. S. A. \textbf{40} (1954), 1147--1151.

\bibitem{BFR86} M. Bordemann, M. Forger and H. Romer; {Homogeneous Kähler
manifolds: Paving the way towards supersymmetric sigma-models}, \emph{%
Commun. Math. Physics},102, (1986), 605--647.

\bibitem{cheeger} J.Cheeger and D.G.Ebin; \emph{Comparsion theorem in
Riemannian geometry, } {North-Holland, Amsterdam-Oxford}; Elsevier, New York
(1975).

\bibitem{CNG} N.Cohen, L.Grama and C.Negreiros; \emph{Equigeodesics on flag
manifolds, }{preprint arXiv:0904.3770} (2009).

\bibitem{ncsa} N.Cohen, C. Negreiros and L. San Martin; \emph{A rank-three
condition for invariant $(1,2)$-symplectic almost Hermitian structures on
flag manifolds, } Bull. Braz. Math. Soc. (N.S.) \textbf{33} (2002), no. 1,
49--73.

\bibitem{neitonnir} N.Cohen and N.P.Silva; \emph{$F$-structure and $T$-roots
on generalized flag manifolds}, in preparation.

\bibitem{MN} X. Mo and C. J. C. Negreiros, $(1,2)$-Symplectic structures on
flag manifolds, Tohoku Math.\ J. \textbf{52} (2000), 271--282.

\bibitem{Evan2} C. J. C. Negreiros and E. C. F. Santos, \textit{Einstein
metrics on flag manifolds}, Revista De La Union on Matemática Argentina,
Vol. 47, Num. 2, 2006.

\bibitem{helg} S. Helgason; \emph{Differential geometry and symmetric
spaces, } Academic Press, New York-London 1962.

\bibitem{H} J. Humphreys, \textit{Introduction to Lie Algebras and
Representation Theory}, Springer, 1972.

\bibitem{kimura} M. Kimura, \textit{Homogeneous Einstein metrics on certain K%
ähler C-spaces}, Adv. Stud. Pure Math. 18-I (1990) 303-- 320.

\bibitem{KN} S.Kobayashi and K.Nomizu; \emph{Foundations of differential
geometry, Vol. II.} Interscience Publishers John Wiley and Sons, Inc., New
York-London-Sydney (1969).

\bibitem{K1} O.Kowalski and J.Szenthe; {On the existence of homogeneous
geodesics in homogeneous riemannian manifolds}, \emph{Geo.Dedicata} \textbf{%
81}(2000), 209-214.

\bibitem{K2} O.Kowalski and L.Vanhecke; {Riemannian manifolds with
homogeneous geodesics}, \emph{Bolletino U.M.I} \textbf{7}(5-B)(1991),
189-246.

\bibitem{sakane} Y. Sakane,Homogeneous Einstein metrics on flag manifolds, 
\emph{Lobachevskii J. Math.}\textbf{4} (1999),71-87.

\bibitem{sm} L.San Martin and C.Negreiros; {Invariant almost Hermitian
structures on flag manifolds,} \emph{Adv. Math.} \textbf{178} (2003),
277--310.

\bibitem{neiton} N. P. da Silva, \textit{Métricas de Einstein e estruturas
Hermitianas invariantes em variedades bandeira}, Tese de Doutorado,
Universidade de Campinas, 2009.

\bibitem{Wa} H. C. Wang, {Closed manifolds with homogeneous complex structure%
}, \emph{Amer. J. Math} \textbf{76} (1954),1-32.

\bibitem{wangziller} M. Wang and W. Ziller, {Existence and non-existence of
homogeneous Einstein metrics},\emph{Invent. Math.} \textbf{84} (1986),
177-194.

\bibitem{Wang e Ziller} M.Wang and W. Ziller, \textit{On normal homogeneous
Einstein metrics}, Ann. Sci. Ecole Norm. Sup. 18 (1985), 563-633.
\end{thebibliography}
\end{document}